\documentclass[a4paper,11pt]{article}
\usepackage{amsmath,amssymb, amsthm}

\input xy
\xyoption{all}

\def\X{\tilde{X}}
\def\d{\tilde{d}}

\def\R{\mathbb{R}}

\def\l{\ell}

\def\b{\bullet}
\def\v{\mathrm{vol}}

\def\p{\underline{p}}
\newtheorem{theo}{\bf{Theorem}}
\newtheorem{lem}[theo]{Lemma}
\newtheorem{cor}[theo]{Corollary}

\newtheorem{prop}[theo]{Proposition}

\begin{document}
\title{Minimal Volume Entropy on Graphs}
\author{Seonhee Lim}
\date{}
\maketitle
\begin{abstract}
 Among the normalized metrics on a graph, we show the existence and the uniqueness of an entropy-minimizing metric, and
give explicit formulas for the minimal volume entropy and the metric realizing it.
\\

Parmi les distances normalis\'ees sur un graphe, nous montrons l'existence et l'unicit\'e d'une distance qui
minimise l'entropie, et nous donnons des formules explicites pour l'entropie volumique minimale et la distance
qui la r\'ealise.
\end{abstract}

\paragraph{1.Introduction} Let $(X,g)$ be a compact connected Riemannian manifold of nonpositive curvature. It was
shown by A. Manning \cite{Man} that the topological entropy $h_{top}(g)$ of the geodesic flow is equal to the
volume entropy $h_{vol}(g)$ of the manifold
$$ h_{vol}(g)= \underset{r \to \infty}{\lim} \frac{1}{r} \log (vol(B(x,r))),$$ where $B(x,r)$ is the ball of radius
$r$ centered at some point $x$ in the universal cover $\X$ of $X$. G. Besson, G. Courtois and S. Gallot
\cite{BCG} proved that if $X$ has dimension at least three and carries a rank one locally symmetric metric
$g_0$, then for every Riemannian metric $g$ such that $vol(X,g)=vol(X,g_0)$, the inequality
$$h_{vol}(g) \geq h_{vol}(g_0)$$
holds, with equality if and only if $g$ is isometric to $g_0$. This solved a conjecture mainly due to M. Gromov
\cite{Gro}, which had been proved earlier by A. Katok \cite{Kat} for metrics in the conformal class of the
hyperbolic metric on a compact orientable surface.

 \vspace{0.1 in} In this paper, we are interested in the analogous
problem for finite graphs, endowed with metrics obtained by varying the length $\l(e)$ of the edges $e \in EX$
of a graph $X$. Regular and biregular trees, as rank one buildings (see \cite{BT}), are non-archimedian analogs
of rank one symmetric spaces, and they carry many lattices (see \cite{BL} for instance). As in the case of
Riemannian manifolds, it is well known that the volume entropy of a finite metric graph is equal to the
topological entropy of the geodesic flow on its universal cover (see \cite{Gui}) and also equal to the critical
exponent of its fundamental group acting on its universal covering tree (\cite{Bou}).

We show that for any graph, among all (normalized) metrics on it, there exists a unique metric minimizing the
volume entropy. We further give explicit formulas for the minimal volume entropy and the metric realizing it.

\medskip
\noindent \textbf{Theorem.} {\it Let $X$ be a finite connected graph such that the valency at each vertex $x$,
which we denote by $k_x +1$, is at least $3$. Then there is a unique normalized length distance minimizing the
volume entropy $h_\v (d)$. The minimal volume entropy is $$h_{\min}= \frac{1}{2}\underset{x \in VX}{\sum}
(k_x+1)\log k_x,$$ and the entropy minimizing length distance $d=d_\l$ is given by
$$\forall e \in EX, \;\;\; \l(e)=\frac{\log (k_{i(e)} k_{t(e)})}{\underset{x \in VX}{\sum} (k_x+1)\log k_x}.$$
}
\medskip

In the special case of a regular graph, which is the analog of a Riemannian manifold carrying a locally
symmetric metric, the volume entropy is minimized by the metric for which all the edges have the same length.
This special case was independently shown by I.~Kapovich and T.~Nagnibeda (\cite{KapNag}).

The above theorem has an analog for graphs of groups (see \cite{Ser}) as well, as in
Proposition~\ref{graphofgroups}. Finally, using Proposition~\ref{graphofgroups}, we show that for a $n$-sheeted
covering graph of groups (see for instance \cite{Bas}) $\phi:~(Y,H_\b) \to (X,G_\b)$, there holds, for the
proper definition of the volume of a metric graph of groups (see section 4),
$$ h_{vol}(Y,H_\b,d)vol(Y,H_\b,d) \geq n \; h_{vol}(X,G_\b,d_0) vol(X,G_\b, d_0),$$ and that the equality holds
 if and only if
the length distance $d$ on $(Y,H_\b)$ is an entropy-minimizing length distance among the length distances of the
same volume, and the map $\phi$ is a metric covering from $(Y,H_\b,d)$ to $(X, G_\b, \lambda d_0),$ for some
$\lambda>0$. This can be considered as an analog of the main theorem in \cite{BCG} for graphs.

\vspace{0.1 in} {\it Acknowledgement.} We are grateful to F. Paulin for suggesting this problem as well as for
helpful discussions and encouragement. We would also like to thank H.~Furstenberg and G.~Margulis for insightful
discussions and H.~Furstenberg for explaining us the relevant work \cite{FK}. We thank I. Kapovich for pointing out the preprint of I. Rivin (\cite{Riv}).

\paragraph{2. Volume entropy and path growth}
Let us consider a nonempty connected \textit{unoriented} finite graph $X$ without any terminal vertex. We will
denote the set of vertices by $V\!X$ and the set of oriented edges of $X$ by $E\!X$. We denote again by $X$ the
geometric realization of $X$. For every edge $e$, let us denote by $i(e)$ and $t(e)$ the initial and the
terminal vertex of $e$, respectively. We define a {\it length distance}  $d$ on $X$ by assigning a positive real
number $\ell (e)=\ell (\bar e)$ for each
 unoriented edge $\{e, \bar e \}$ of $X$, and by letting $d=d_\ell :X \times X \to [0, \infty [ $ be the maximal
 distance which makes each half-edge of an edge $e$ containing a vertex, isometric to $[0, \frac{\ell(e)}{2}]$.
For a length distance $d_\l$, let $l_{\max}=\underset{e \in EX}{\max}\; \l(e)$ and $l_{\min}=\underset{e \in
EX}{\min}\; \l(e)$. Define the volume of $X$ by $$vol(X,d)=\frac{1}{2} \underset{e \in EX}{\sum} \l(e),$$ i.e.,
the sum of lengths of the {\it unoriented} edges. We denote by $\Delta(X)$ the set of all length distances
$d=d_\l$ on $X$ normalized so that $vol(X,d)=1.$

 For a fixed length
distance $d$, let us consider a universal covering tree $\X \to X$ equipped with the lifted distance $\tilde d$
of $d$. For any connected subset $S$ of $\tilde X$, let us denote by $\ell(S)$ the sum of the lengths of (the
maximal pieces of) the edges in $S$. We define the {\it volume entropy} $h_{vol}(d)=h_\v (X, d)$ as
$$ h_{vol}(d)=\underset {r \to \infty}{ \limsup}\; \frac{1}{r} \log \ell(B(x_0, r)),$$
where $B(x_0, r)=B_d(x_0,r)$ is the ball of radius $r$ with center a fixed vertex $x_0$ in $(\X, \d)$. The
entropy $h_{vol} (d)$ does not depend on the base point $x_0$, and we may sum either on the oriented or on the
non-oriented edges. Note also the homogeneity property $$ h_{vol} (d_{\alpha \ell}) = \frac{1}{\alpha} h_{vol}
(d_\ell), \;\;\;\;\;\;\; (*)$$ for every $\alpha>0$. Remark that $h_\v (X,d)\v(X,d)$ is invariant under
dilations, therefore to minimize the entropy with constant volume, it suffices to consider the length metrics of
volume $1$.

If $\pi_1X$ is not cyclic, or equivalently if $X$ has no terminal vertices and is not reduced to one cycle, then
$h_{\v}=h_\v (d)$ is strictly positive, which we will assume from now on (see for instance \cite{Bou}). It was
shown by Roblin (\cite{Robl}) that the upper limit above is in fact a limit. This implies that as $r \to
\infty$, $$\l(B(x_0, r))=e^{h_{\v(d)} r+o(r)}.$$

By a {\it metric path of length $r$ in $X$}, we mean the image of a local isometry $f:[0,r] \to X$. Note that
the endpoint of a metric path is not necessarily a vertex. By a {\it combinatorial $n$-path of length $r$ in
$X$}, we mean a path $\underline{p}=e_{1}e_{2}\cdots e_{n}$ of consecutive edges in $X$ without backtracking
such that $\sum_{j=1}^{n-1} \ell(e_{j})< r \leq \sum_{j=1}^{n} \ell(e_{j})$. A combinatorial path is a
combinatorial $n$-path for some $n$.

\begin{lem}\label{first} Let $N_r(x_0)$ be the cardinality of the set of combinatorial paths of length
$r$ in $\X$ starting at $x_0 \in V\X$. Then the number $N_r(x_0)$ satisfies
$$\underset{r \to \infty}{\limsup}\frac{\log N_r (x_0)}{r}=\underset{r \to \infty}{\lim}\frac{\log N_r (x_0)}{r}=h_{vol}.$$
\end{lem}

\noindent {\it Proof.} It follows directly from $\l(B(x_0, r))=e^{(h_\v+o(1))r}$ that for any $l >0$,
$$\underset{r \to \infty}{\limsup}\frac{\log \l(B(x_0,r)-B(x_0,r-l))}{r}=\underset{r \to \infty}{\lim}\frac{\log \l(B(x_0,r)-B(x_0,r-l))}{r}=h_{vol}.$$
Now let $N'_r(x_0)$ be the cardinality of the set of metric paths of length $r$ starting at $x_0$. As $\X$ has
no terminal vertices, for any $l>0$,
$$ l N'_{r-l} (x_0) \leq \l(B(x_0,r)-B(x_0,r-l)) \leq l N'_r (x_0).$$
Therefore
$$\underset{r \to \infty}{\limsup}\frac{\log N'_r(x_0)}{r}=\underset{r \to \infty}{\lim}\frac{\log N'_r(x_0)}{r}=h_{vol}.$$
It is clear that we get a combinatorial path of length $r$ by continuing a metric path of length $r$ until it
meets a vertex. Also, two distinct combinatorial paths of length $r$ cannot be extensions of one metric path of
length $r$ by the strict inequality in the definition of a combinatorial path. It follows that $ N_r(x_0) = N'_r
(x_0),$ thus $N_r(x_0)$ has the same exponential growth rate as $N'_r(x_0)$, which is $h_{vol}$. \qed

\vspace{0.1 in} \noindent Let $A=A(X)$ be the {\it edge adjacency matrix of $X$}, i.e. a $|EX| \times |EX|$
matrix defined by $A_{ef} =\rho_{ef}$, where $\rho_{ef}$ has value $1$ if $ef$ is a combinatorial path, i.e. if
$t(e)=i(f)$ and $\bar{e} \neq f$, and value 0 otherwise. It is easy to see that the entry $A^n _{ef}$ is nonzero
if and only if there is a combinatorial $n$-path starting with $e$ and ending with $f$. (Note that the
definition of $\rho_{ef}$ implies that such a path does not have backtracking.)

Let us show that for any connected graph without any terminal vertex, which is not a cycle, the matrix $A$ is
irreducible. Recall that a matrix $M$ is {\it reducible} if there exists a permutation matrix $P$ such that
$PMP^{-1}$ is a block diagonal matrix of at least two nontrivial blocks. A matrix is {\it irreducible} if it is
not reducible. It is clear that a nonnegative matrix $M$ is reducible if and only if $M+M^t$ is reducible. The
matrix $A +A^t$ has entries $b_{ef}$ with $b_{ef}$ nonzero if either $ef$ or $fe$ is a combinatorial path. Let
$\sim$ be the equivalence relation on $EX$ generated by the relation $e \sim f$ if $\rho_{ef}=1$ (i.e if $ef$ is
a path without backtracking).

\begin{lem} If $EX/\sim$ consists of only one element, then $A$ is irreducible.
\end{lem}

\noindent {\it Proof.} By contradiction. Assume that $A$ is reducible, and hence $A+A^t$ is also reducible. This
means that it is possible to find two complementary subsets $U$ and $V$ in $EX$ such that $b_{uv}=0$ for any $u
\in U$ and $v \in V$. But then clearly no element $u$ of $U$ can be equivalent to an element $v$ of $V$, and in
particular there are at least two equivalence classes for the equivalence relation $\sim$. \qed

\begin{prop} Let $X$ be a connected graph without any terminal vertex. Then the matrix $A$ is irreducible if and
only if $X$ has a vertex of valency at least three.
\end{prop}

\noindent \textit{Proof.} We first claim that for any $e,f \in EX$, either $e \sim f$ or $e \sim \overline{f}$.
Let $\mathbf{e}$ and $\mathbf{f}$ be the unoriented edges underlying $e$ and $f$. There is a shortest (thus
without backtracking) unoriented path in $X$ linking $\mathbf{e}$ and $\mathbf{f}$. Thus, there is an
combinatorial (oriented) path in $X$ linking either $e$ or $\overline{e}$ to either $f$ or $\overline{f}$. If
$\underline{p}=e_1 \cdots e_s$ is a combinatorial path, then so is $\overline{\underline{p}}=\overline{e_s}
\cdots \overline{e_1}$. We deduce that, for any two edges $e_1$ and $e_2$, $e_1 \sim e_2 \Leftrightarrow
\overline{e_2} \sim \overline{e_1}$. Applying this to $e, \overline{e}$ and $f, \overline{f}$ yields the claim.

Now if $X$ has only bivalent vertices (i.e. if $X$ is a cycle, as $X$ has no terminal vertices) then it is easy
to see that $A$ is not irreducible. Conversely, let $X$ be connected and let $x$ be a vertex of valence at least
three, with outgoing edges $e_1, e_2$ and $e_3$ (these may be loops). Then $\overline{e_1} \sim e_3$ and
$\overline{e_2} \sim e_3$ hence $\overline{e}_1 \sim \overline{e_2}$. But $\overline{e_1} \sim e_2$, and
therefore $e_2 \sim \overline{e_2}$. Now for any edge $e$ in $EX$, either $e \sim e_2$ or $e \sim
\overline{e_2}$ by the claim above. In both cases, we have $e \sim e_2$. This implies that $EX/\sim$ has only
one element and $A$ is irreducible. \qed

\vspace{0.1 in} \noindent Now consider the matrix $A'=A'(d,h)$ defined by $A'_{ef}=\rho_{ef} e^{-h\l(f)}$,
depending on $h$ and the length distance $d_\l$ on $X$. The matrix $A'$ is clearly irreducible since $A$ is
irreducible.

\begin{theo}\label{pathgrowth} Let $X$ be a connected finite graph without any terminal vertex, which is not a
cycle, endowed with a length distance $d=d_\l$. The volume entropy $h_{\v}$ is the only positive constant $h$
such that the following system of linear equations with unknowns $(x_e)_{e \in EX}$ has a solution with $x_e>0$
for every $e \in EX$.
$$ x_e = \sum_{f \in EX} \rho_{ef} e^{-h\ell(f)}x_f, \;{\it for\;all\;} e\in EX. \;\;\;\;\;(**)$$
\end{theo}

\noindent {\it Proof.} By the assumption on the graph, for every $h>0$, we can apply Perron-Frobenius theorem
(see [Gan] for example) to the irreducible nonnegative matrix $A'=(\rho_{ef} e^{-h\ell(f)})$, which says that
the spectral radius of the matrix $A'(h)$ is a positive eigenvalue $\lambda(h)$, which is simple, with an
eigenvector $(x_e=x_e(h))$ whose entries are all positive. The function $\lambda : \R_{\geq 0} \to \R_{\geq 0}$
is clearly a continuous function of $h$ since the characteristic function of the matrix $A'$ is a polynomial of
$\{ e^{-h\l(e)} : e \in EX \}$, and $\lambda(0) \geq 1$ since $\lambda(0)$ is the spectral radius of an
irreducible nonzero matrix $A'(0)$ of nonnegative integer coefficients. Also, $\lambda(h) \to 0$ as $h \to
\infty$, since the coefficients of $A'(h)$ tends to $0$ as $h \to \infty$. By the mean value theorem, there
exists an $h$ satisfying $\lambda(h)=1$.

 Now assume that $h>0$ satisfies $(**)$ for some positive $x_e$'s.
 Fix an arbitrary edge $e \in EX$ (any edge $e$ satisfies $x_e > 0$ by Perron-Frobenius theorem),
and choose a vertex $x_0$ in $\X$ which is an initial vertex of a fixed lift $\tilde{e}$ of $e$ in $\X$. Let us
fix a positive constant $r$.

 Let $P_r(e)$ be the set of combinatorial paths of length $r$ in $X$ starting with $e$. We will denote a
combinatorial path in $X$ by $\p=e_1 e_2 \cdots e_n$, its terminal edge by $t(\p)=e_n$ and its metric length by
$\l(\underline{p})= \sum_{i=1}^{n} \l(e_i)$. Denote by $\mathcal{P}_n(e)$ (resp. $\mathcal{P}'_n(e)$) the set of
combinatorial $k$-paths of length $r$ with $k \leq n$ (resp. combinatorial $n$-paths of length strictly less
than $r$) in $X$ starting with $e$. Remark that $\mathcal{P}_n(e) \cap \mathcal{P}'_n(e)=\emptyset$ and if $n$
is large enough, $\mathcal{P}_n(e)=P_r(e)$ and $\mathcal{P}'_n(e)=\emptyset$.

Let us rewrite the equation (**) as
$$ e^{h\ell(e)}x_e=\underset{ \p \in \mathcal{P}_{2}(e) \cup \mathcal{P}'_2(e)}{\sum} e^{-h\ell(\p)}x_{t(\p)}.$$

 Let us replace each $x_{t(\p)}$ in the above equation by $\underset{f
\in EX}{\sum} \rho_{t(\p)f}e^{-h\ell(f)}x_f$ whenever $\ell(\p) < r$, i.e. when $\p \in \mathcal{P}'_2(e)$. The
resulting equation is
$$ e^{h\ell(e)}x_e=\underset{ \p \in \mathcal{P}_{3}(e) \cup \mathcal{P}'_3(e)}{\sum} e^{-h\ell(\p)}x_{t(\p)}.$$

Repeat this process: at each step, for each $\underline{p} \in \mathcal{P}'_n(e)$, replace $x_{t(\p)}$ on the
right hand side of the previous equation by $\underset{f \in EX}{\sum} \rho_{t(\p)f}e^{-h\ell(f)}x_f$, to get
$$e^{h\l(e)}x_e=\underset{\p \in \mathcal{P}_{n+1}(e) \cup \mathcal{P}'_{n+1}(e)}{\sum} e^{-h\l(\p) x_{t(\p)}}.$$ For $n$ large enough, the
resulting equation is
$$e^{h\l(e)} x_e  = \underset{\p \in P_r(e)}{\sum}  e^{-h\l(\p)} x_{t(\p)}.$$

(In the case when the lengths of the edges are all equal to $1$ and $r$ is a positive integer, we continue until
we get the equation $\underline{x}= A^{r-1} \underline{x}$.)

 Then in the resulting equation, the number of times each $x_f$ appears on the right hand side is exactly the
number $N_r(e,f)$ of combinatorial paths of length $r$ in $\X$ with initial edge $\tilde{e}$ and terminal edge
some lift of $f$ in $\X$. Note also that the metric length of such a path is at least $r$ and less than
$r+l_{\max}$. Thus

$$ \underset{f \in EX}{\sum} N_r(e,f) e^{-h(r+l_{\max})} x_f \leq e^{h\l(e)} x_e \leq  \underset{f \in EX}{\sum} N_r(e,f)
e^{-hr} x_f. $$

Now if $h$ is strictly greater than the volume entropy $h_{\v}$, then
\begin{align*}
0 < e^{h\l(e)}x_e &\leq \sum N_r(e,f) e^{-hr}x_f \leq N_r(x_0) e^{-hr} \sum x_f \\
& \leq e^{(h_{\v}+0(1))r} e^{-hr} \sum x_f = e^{r(h_{\v}-h+o(1))} \sum x_f \to 0
\end{align*}
as $r$ goes to infinity, which is a contradiction. On the other hand, suppose that $h$ is strictly smaller than
the volume entropy $h_{\v}$. As $N_r(x_0)= \underset{{e,f \in EX,}{\;i(\tilde{e})=x_0}}{\sum} N_r(e,f)$, there
exist some $e$ and $f$, depending on $r$, such that $ N_r(e,f) e^{-hr} \geq \frac{1}{|EX|^2} N_r(x_0) e^{-hr}$.
Since
\begin{align*}
e^{h\l(e)} x_e &\geq  N_{r}(e,f) e^{-h(r+l_{\max})} x_f \geq \frac{1}{|EX|^2} e^{(h_\v+o(1))r}
e^{-h(r+l_{\max})}
x_f\\
&\geq \frac{1}{|EX|^2}e^{(h_{\v}-h+o(1))r} e^{-hl_{\max}} \underset{f \in EX}{\min}\{x_f\},
\end{align*}
 and
 the latter goes to infinity as $r$ goes to infinity, it follows that $x_e=\infty$, which is again a contradiction.
We conclude that $h$ is equal to the volume entropy $h_{\v}$. \qed

\vspace{0.1 in} \noindent {\it Remark.} Hersonsky and Hubbard showed in \cite{HerHub} that the Hausdorff
dimension of the limit set of a Schottky subgroup of the automorphism group of a simplicial tree satisfies
similar systems of equations.

\paragraph{3. Minimal volume entropy} In this section, we prove the main theorem announced in the introduction, using
Theorem~\ref{pathgrowth}.
\begin{theo}\label{entropy} Let $X$ be a finite connected graph such that the valency at each
vertex $x$, which we denote by $k_x +1$, is at least $3$. Then there is a unique $d$ in $\Delta(X)$ minimizing
the volume entropy $h_\v (d)$. The minimal volume entropy is $$h_{\min}(X)= \frac{1}{2}\underset{x \in VX}{\sum}
(k_x+1)\log k_x,$$ and the entropy minimizing length distance $d=d_\l$ is characterized by
$$\l(e)=\frac{\log k_{i(e)}k_{t(e)}}{\underset{x \in VX}{\sum} (k_x+1)\log k_x}, \;\;\;\forall e \in EX.$$
\end{theo}

\noindent {\it Remark.} Since we can eliminate all the vertices of valency two without changing the entropy, the
existence of $d$ in $\Delta(X)$ minimizing the volume entropy, with minimal value given by the same formula,
holds for any graph who does not have a terminal vertex and is not isometric to a circle. What is uniquely
defined at such a minimum is the length of each connected component of $X$ where the vertices of valency at
least three are removed.

\vspace{.1 in} \noindent {\it Proof.} By assumption, $k_x \geq 2$ for every $x \in VX$. By
Theorem~\ref{entropy}, the volume entropy $h=h_\v$ satisfies
\begin{equation}\label{inequal}
 x_{e} = \underset{f \in EX}{\sum} \rho_{ef} e^{-h\ell(f)} x_f,
\end{equation}
  for each edge $e \in EX$ for some positive
$x_e$'s. Set $y_{e}= e^{-h \ell(e)}x_{e}>0$ for each edge $e$. Then the above equations implies
$$ e^{h\ell(e)} y_{e} = \sum_{f \in EX}\rho_{ef} y_{f} \geq k_{t(e)} \underset{{f \in EX,}\;{\rho_{ef}=1}}{\prod} y_{f}^{1/k_{i(f)}}.$$
The last inequality is simply the inequality between the arithmetic mean and the geometric mean of $y_{f}$'s,
since there are exactly $k_{t(e)}=k_{i(f)}$ edges $f$ such that $\rho_{ef}=1$. Multiplying over all the edges,
we get
$$
 \underset{e \in EX}{\prod} e^{h\ell_{e}} y_{e} \geq \underset{e \in EX}{\prod}( k_{t(e)} \underset{{f \in
EX,}\;{\rho_{ef}=1}}{\prod} y_{f}^{1/k_{i(f)}}).
$$

On the right hand side of the equation, each term $y_f ^{1/k_{i(f)}}$ appears exactly $k_{i(f)}$ times, since
each edge $f$ follows exactly $k_{i(f)}$ edges with terminal vertex $i(f)$. Canceling $\underset{e\in EX}{\prod}
y_{e}>~0$ from each side, we get

\begin{equation}
\label{second} e^{2h} \geq \underset{e \in EX}{\prod} k_{t(e)}=\underset{x \in VX}{\prod} k_x ^{(k_x+1)},
\end{equation}
 since $\underset{e \in EX}{\sum} \ell(e)=2$. The equality holds if and only if equality in the inequality
 (\ref{inequal})
 holds for each $e \in EX$, i.e. the $y_{f}$'s, for $f \in EX$ following $e$, are all equal.

Suppose that the equality in the inequality (\ref{second}) holds. In particular, $$h=\frac{1}{2}\underset{x \in
VX}{\sum}(k_x +1)\log k_x.$$ Since the valency at each vertex is at least $3$, we can choose another edge $g
\neq f$ followed by $e$ and conclude that $y_{f}$ depends only on the initial vertex $i(f)$ of $f$. Let
$z_{i(f)}=y_f>0$. Then the equation $(**)$ in Theorem~\ref{pathgrowth} amounts to
$$e^{h\l(e)}z_{i(e)}=\underset{f \in EX}{\sum} \rho_{ef}z_{i(f)}=k_{t(e)}{z_{t(e)}}.$$ Since $\l(e)=\l(\bar{e})$,
we also have $e^{h\l(e)}z_{t(e)}=k_{i(e)}z_{i(e)}$. Thus $z_{i(e)}/z_{t(e)}=k_{t(e)}/e^{h\l(e)}=
e^{h\l(e)}/k_{i(e)}$ and $$e^{h \l(e)}=\sqrt {k_{i(e)}k_{t(e)}},$$ so that

\begin{equation}\label{third}
\l(e)=\frac{\log k_{i(e)}k_{t(e)}}{ \underset{x \in VX}{\sum}
 (k_x +1)\log k_x}.
\end{equation}
In particular, $\l$ is uniquely defined by this formula. The length distance defined by the formula
(\ref{third}) clearly satisfies the equations $(**)$,
 with
$$h= \frac{1}{2}\underset{x \in VX}{\sum} (k_x+1)\log k_x,$$ and $x_e$'s defined, uniquely up to constant, by setting
$$\frac{e^{-h\l(e)}x_e}{e^{-h\l(f)}x_f }=\sqrt{\frac{k_{t(e)}}{k_{i(e)}}},$$
 for every $f$ such that $i(f)=t(e)$. It is clearly well-defined, since if there is a cycle consisting of
consecutive edges $(e_1, e_2, \cdots, e_n, e_{n+1}=e_1)$, then
\begin{align*}
y_{e_n} = y_{e_{n-1}}\sqrt{\frac{k_{i(e_{n-1})}}{k_{i(e_{n})}}} =\cdots = y_{e_{1}}\prod_{j=2}^{j=n}
\sqrt{\frac{k_{i(e_{j-1})}}{k_{i(e_j)}}} =y_{e_1} \sqrt{\frac{k_{i(e_1)}}{k_{i(e_{n})}}}.
\end{align*}

 By uniqueness in
Theorem~\ref{pathgrowth}, the positive number $h$ given above is the volume entropy of the given length
distance, and it is the minimal entropy of the graph.\qed
\begin{cor}\label{biregular}
If $X$ is a $(k_1+1,k_2+1)$-biregular graph, with $k_1 >1, k_2 >1$, then the volume entropy of the normalized
length distances on $X$ is minimized exactly when the lengths of the edges are all equal, and the minimal volume
entropy is $\frac{|EX|}{4}\log (k_1 k_2)$.
\end{cor}

\noindent {\it Proof.} Suppose that $X$ a $(k_1+1, k_2+1)$-biregular graph, i.e. $k_{i(e)}k_{t(e)}=k_1 k_2$ for
any edge $e$. Let $d=d_\l \in \Delta(X)$ be the entropy-minimizing length distance. Then $\l(e)=\frac{1}{2h}
\log (k_1 k_2)$ does not depend on $e$, thus $\l(e)=\frac{2}{|EX|}$. From $e^{h \l(e)}=\sqrt
{k_{i(e)}k_{t(e)}},$ the volume entropy of this length distance is $h=\frac{|EX|}{4} \log (k_1 k_2)$. \qed

\begin{cor}\label{regular}
If $X$ is a $(k+1)$-regular graph, with $k >1$, then the volume entropy of the normalized length distances on
$X$ is minimized exactly when the lengths of the edges are all equal, and the minimal volume entropy is
$\frac{|EX|}{2}\log k$.
\end{cor}

\noindent {\it Proof.} This is a special case of the above corollary with $k_1=k_2=k$.\qed

\vspace{0.1 in} \noindent {\it Remark.} The last corollary appears implicitly in a preprint of I. Rivin
(\cite{Riv}). There he considers graphs with weights given on the vertices rather than the edges. The directed
line graph $L(X)$ of a graph $X$ is an oriented graph defined so that $VL(X)=EX$ and $EL(X)=\{ (a,b) \in EX^2 :
t(a)=i(b), a \neq \bar{b} \}$. To a given set of weights on the edges $\{\l(e)\}_{EX}$, is associated a set of
weights $\{\l'(x)\}_{VL(X)}$ on the vertices of $L(X)$. One can sees that paths on $X$ without backtracking
correspond to paths with backtracking on $L(X)$, see \cite{Riv} page 14. The minimum of volume entropy of the
graph $L(X)$ with vertex weights $h((\l'(x)))_{VL(X)}$ (computed by I. Rivin) lies in the image of the map
$(\l(e)) \mapsto (\l'(x))$ only when the graph is regular. It seems that for general graphs, one result cannot
be deduced from the other.

\vspace{0.1 in} \noindent {\it Remark.} Corollary~\ref{regular} was also shown independently by I. Kapovich and
T. Nagnibeda \cite{KapNag} by a different method (using random walks). Note that one of their main results, on
the minimal entropy among all graphs having a fixed fundamental group, can be deduced from Theorem~\ref{entropy}
as in the following corollary. A special case when the graph has a highly transitive automorphism group had been
shown earlier by G. Robert (\cite{Rob}).

\begin{cor}(\cite{KapNag} Theorem B) Consider the set of all finite metric graphs without a vertex of valency
one or two, whose fundamental group is a free group of given rank $r\geq 2$. Then among volume 1 length metrics,
the volume entropy is minimized by any (regular) trivalent graph in this set, with the metric assigning the same
length for every edge.
\end{cor}

\noindent {\it Proof.} Let $(X,d)$ be such a graph. Suppose that there is a vertex $x$ of valency $k_x+1$
strictly greater than three, with outgoing edges $e_1, \dots, e_{k_x+1}$. Let us introduce a new vertex $y$ and
a new edge $f$, and replace $x$ and its outgoing edges $e_1, \cdots, e_{k_x +1}$, by two vertices $x$ and $y$,
with outgoing edges $f, e_3, \cdots, e_{k_x +1}$ and $e_1, e_2, \bar{f}$, respectively. Repeat the operation on
$x$, until the valency of $x$ reduces to three, to get a new graph $X'$. The graph $X'$ has $k_x - 2$ more
vertices than $X$, all with valency three.

Let $d_0$ and $d_0'$ be the unique normalized entropy-minimizing length distances on $X$ and $X'$, respectively.
By the formula in Theorem~\ref{entropy}, since for $t\geq 3$, $(t+1)\log t> (t-1)3 \log 2$, it follows that
\begin{align*}
h_{\v}(X,d) &\geq h_{\v}(X,d_0)= \frac{1}{2}\underset{z \in VX-\{x\}}{\sum} (k_z+1)\log k_z+(k_x+1)\log k_x \\
&> \frac{1}{2}\underset{z \in VX-\{x\}}{\sum} (k_z+1)\log k_z+(k_x-1)3\log2=h_{\v} (X',d'_0).
\end{align*}
Repeat the operation until we get a regular trivalent graph. Now by Corollary~\ref{regular}, the volume entropy
is minimized when all the edges have the same length. \qed

\vspace{.1 in}

\paragraph{4. Entropy for graphs of groups}
\noindent As another corollary of Theorem~\ref{entropy}, let us show the analogous result of
Theorem~\ref{entropy} for graphs of groups. Let $(X,G_\b)$ be any finite connected graph of finite groups.
(Basic references for graphs of groups are \cite{Ser} and \cite{Bas}.) Let $T$ be a (Bass-Serre) universal
covering tree of $(X,G_\b)$ and let $p: T \to X$ be the canonical projection. The \textit{degree of a vertex $x$
of $(X,G_\b)$} is defined by
$$\sum_{e \in EX, i(e)=x}\frac{|G_{x}|}{|G_e|}.$$ It is easy to see that it is equal to the valency of any lift of $x$
in $VT$, and we will denote it again by $k_x +1$. We define a \textit{length distance $d_{\l}$ on $(X, G_\b)$}
as a length distance $d_\l$ on the underlying graph $X$. The \textit{volume of $(X,G_\b,d_\l)$} for a given
length distance $d_\l$ on $(X,G_\b)$, is defined by
$$vol_\l(X,G_\b)= \frac{1}{2} \sum_{e \in EX} \frac{\l(e)}{|G_e|}.$$ Note that in the
case where $\l(e)$ is equal to $1$ for every edge $e$ and $T$ is $k$-regular, the volume $vol_\l(X,G_\b)$ is
$k/2$ times the usual definition of the volume $\sum_{x \in VX} 1/|G_x|$ of a graph of groups since
$k=\sum_{e\in EX, i(e)=x}|G_x|/|G_e|$. The \textit{volume entropy $h_\v (X,G_\b ,d_\l )$ of $(X,G_\b,d_\l)$} is
defined to be the exponential growth of the balls in $T$ for the lifted metric as in the case of graphs.

\begin{prop}\label{graphofgroups} Let $(X,G_\b)$ be a finite graph of finite groups such that the degree at each vertex $x$
 of $(X,G_\b)$ is at least three. Among the normalized (i.e. volume one) length distances on $(X,G_\b)$, there exists a unique
 normalized length distance minimizing the volume entropy. At this minimum, the length of each edge is proportional to
$ \log( k_{i(e)} k_{t(e)})$ and the minimal volume entropy is $$h_{\min}(X,G_\b)=\frac{1}{2} \underset{x \in
VX}{\sum}\frac{(k_x+1) \log k_x}{|G_x|}.$$
\end{prop}
\textit{Proof.} Let $\Gamma$ be a fundamental group of the graph of groups $(X,G_\b)$. There exists a free
normal subgroup $\Gamma'$ of $\Gamma$ of finite index (see \cite{Ser}), say $m$. The group $\Gamma'$ acts freely
on $T$, hence the quotient graph $X'=\Gamma ' \backslash T$ is a finite connected graph. It is easy to see that
each $x$ in $VX$ (resp. $e$ in $EX$) has $\frac{m}{|G_x|}$(resp. $\frac{m}{|G_e|}$) lifts in $VX'$ (resp. $EX'$)
by the canonical map $\pi:X' \to X$, since
$$m=[\Gamma: \Gamma']=\frac{\sum_{x' \in VX'}1}{\sum_{x \in VX} 1/|G_x|}$$ (see \cite{Bas} for example).
 It is clear that $\l'(e)=\l(\pi(e))$ and the valency $k_y +1$ is equal to the
degree $k_{\pi(y)}+1$. Any length distance $d_\l$ of volume one on $(X,G_\b)$ can be lifted to $X'$ to define a
length distance $d_{\l'}$ normalized so that
$$vol_{\l'}(X') = \frac{1}{2} \sum_{e \in EX'}\l'(e)= \frac{1}{2} \sum_{e \in EX} \frac{m}{|G_e|}\l(e)=m.$$ The
volume entropy of $(X', d'_\l)$ is equal to the volume entropy of $(X,G_\b,d_\l)$ as they have the same
universal covering metric tree.  By the homogeneity property $(*)$, we can apply Theorem~\ref{entropy} to
conclude that among the length distances of volume $m$ on $X'$, there exists a unique entropy-minimizing length
distance $d'_0=d_{\l'}$ on $X'$. By uniqueness in Theorem~\ref{entropy}, the length distance $d'_0$ is invariant
under the group $\Gamma / \Gamma '$. In particular, there is a normalized length distance $d_0=d_\l$ on
$(X,G_\b)$ whose lift to $X'$ defines $d'_0$. The minimal volume entropy of $(X,G_\b)$ is clearly the volume
entropy of $(X',d'_0)$ since for any length distance $d$ on $(X,G_\b)$,
$$ h_\v (X,G_\b, d) =h(X',d') \geq h(X',d'_0)=h_\v (X,G_\b,d_0),$$
where $d'$ is defined by the lift of $d$ on $X'$. Since the length $\l'(e)$ of an edge $e$ is proportional to
$\log(k_{i(e)} k_{t(e)})=\log(k_{\pi(i(e))}k_{\pi(t(e))})$ for every edge $e$ in $EX'$, so is true for every
edge $e$ in $EX$. Since each vertex $x$ in $VX$ appears $\frac{m}{|G_x|}$ times in $X'$ and the degree $k_x+1$
is equal to the valency $k_{x'}+1$ of any lift $x' \in \pi^{-1}(x)$ of $x$ in $X'$, the minimal volume entropy
of $(X,G_\b)$ is
\begin{align*}
h_{d_0}(X,G_\b) &=h(X',d'_0)=\frac{1}{m} h(X', \frac{1}{m} d'_0)
               =\frac{1}{2m} \underset{x' \in VX'}{\sum} (k_x'+1) \log k_x'\\
               & =\frac{1}{2m} \underset{x \in VX}{\sum} \frac{m}{|G_x|} (k_x+1) \log k_x
=\frac{1}{2} \underset{x \in VX}{\sum}\frac{(k_x+1) \log k_x}{|G_x|}.
\end{align*}

\qed

\vspace{.1 in} \noindent Now we want to consider a more general situation than in
Proposition~\ref{graphofgroups}. The main theorem in \cite{BCG} says that if $f:(Y,g) \to (X,g_0)$ is a
continuous map of non-zero degree between compact connected $n$-dimensional Riemannian manifolds and $g_0$ is a
locally symmetric metric with negative curvature, then
$$h^n(Y,g)\v(Y,g) \geq |\mathrm{deg}\; f| h^n(X,g_0) \v(X,g_0),$$ and the equality holds if and only if $f$ is homotopic to
a Riemannian covering.

 Let $(X,G_\b, d_0=d_\l)$ be a finite (connected) graph of finite groups endowed with the
normalized length distance minimizing the volume entropy. Let $(Y,H_\b, d)$ be another finite graph of finite
groups with a length distance. Let $\phi=(\phi, \phi_\b, \gamma_\b):~(Y,H_\b)~\to~(X, G_\b)$ be a (Bass-Serre)
covering of graphs of groups (see \cite{Bas}). The value $$n: =\underset{y \in \phi^{-1}(x)}
 \sum \frac{|G_x|}{|H_y|} = \underset{f \in \phi^{-1}(e)} \sum \frac{|G_e|}{|H_f|}$$ does not depend on the vertex $x$
 nor on the edge $e$ of $X$ since the graph $X$ is connected, and it is an integer. A covering graph of groups with the above $n$ is said
 to be {\it $n$-sheeted} (see \cite{Lim}).

When $(Y,H_\b)$ and $(X, G_\b)$ are graphs (of trivial groups), the next corollary can be considered as an
analog of the main theorem in \cite{BCG}.

\begin{cor}
Let $\phi:(Y,H_\b) \to (X,G_\b)$ be a $n$-sheeted covering of graphs of groups and let $d_0$ be the
entropy-minimizing length distance on $(X,G_\b)$ of volume one. Suppose that the degree at each vertex of
$(X,G_\b)$ and $(Y,H_\b)$ is at least three. Then there holds
$$ h_{vol}(Y,H_\b,d)vol(Y,H_\b,d) \geq n \; h_{vol}(X,G_\b,d_0) vol(X,G_\b, d_0).$$ The equality holds if and only if
the length distance $d$ on $(Y,H_\b)$ is a length distance minimizing entropy among the length distances of the
same volume, and in that case the map $\phi$ is a metric covering from $(Y,H_\b,d)$ to $(X, G_\b, \lambda d_0),$
for some $\lambda
>0$.
\end{cor}
\noindent \textit{Proof.} By the homogeneity property $(*)$, we may assume that $\v(Y,H_\b, d)=1$. Applying
Proposition~\ref{graphofgroups} to $(Y,H_\b)$ and $(X,G_\b)$, it follows that there exists a unique length
distance $d'_0=d_{\l'}$ on $Y$ minimizing the volume entropy and that
\begin{align*}
h_\v(Y,H_\b,d) &\geq h_{\min}(Y, H_\b) =\frac{1}{2} \underset{y \in VY}{\sum}\frac{(k_y+1) \log k_y}{|H_y|}
=\frac{1}{2}\underset{x \in VX}{\sum} \underset{y \in \phi^{-1}(x)}{\sum}\frac{(k_x+1) \log k_x}{|H_y|}\\
&=\frac{1}{2}n \underset{x \in VX}{\sum}\frac{(k_x+1) \log k_x}{|G_x|}=nh_{\min}(X,G_\b)=nh_\v(X,G_\b, d_0).
\end{align*}
By Proposition~\ref{graphofgroups}, the equality holds if and only if $d=d'_0$. In that case, the length of each
edge $e$ in $EY$ is proportional to $\log(k_{i(e)}k_{t(e)})=\log(k_{i(\phi(e))}k_{t(\phi(e))})$, thus
proportional to the length of the edge $\phi(e)$. More precisely, let $\l'(e)=c' \log(k_{i(e)}k_{t(e)})$ for
every $e \in EY$ and let $\l(e)=c \log(k_{i(e)}k_{t(e)})$ for every $e \in EX$. From the assumption
$\v_\l(X,G_\b)=\v_{\l'}(Y,H_\b)=1$, it follows that
\begin{align*}
1 &= \frac{1}{2} \underset{ g \in EY}{\sum} \frac{c' \log(k_{i(g)}k_{t(g)})}{|H_g|}
 =\frac{1}{2}\underset{e \in EX}{\sum} \underset{ g \in \phi^{-1}(e)}{\sum}\frac{c' \log(k_{i(g)}k_{t(g)})}{|H_g|}
= \frac{1}{2} \underset{e \in EX}{\sum}  \frac{nc'\log(k_{i(e)}k_{t(e)})}{|G_e|},
\end{align*}
and therefore $$c'=\frac{1}{\frac{n}{2} \underset{e \in EX}{\sum}
\frac{\log(k_{i(e)}k_{t(e)})}{|G_e|}}=\frac{c}{n},$$ in other words, $\l'(e)= \l(e)/n.$

We conclude that for any length distance $d$ on $(Y,H_\b)$, there holds $$ h_{vol}(Y,H_\b,d)vol(Y,H_\b,d) \geq n
\; h_{vol}(X,G_\b) vol(X,G_\b, d_0).$$ By Proposition~\ref{graphofgroups} the equality holds if and only if $d$
is proportional to $d'_0$, say $d=\lambda n d'_0$ for some $\lambda>0$. Then the length of each edge $e$ in
$(Y,H_\b,d)$ is $\lambda \l(\phi(e))$, and the map $\phi$ is a metric covering from $(Y,d'_0)$ to $(X, \lambda\;
d_0) $. \qed

 \small{

\paragraph{}

\quad Yale University,  New Haven, CT 06520-8283, USA and

ENS-Paris,  UMR 8553 CNRS, 45 rue d'Ulm, 75230  Paris Cedex 05,  France

seonhee.lim@yale.edu, Seonhee.Lim@ens.fr

\end{document}